\documentclass{amsart}
\usepackage{pstricks}
\usepackage{amsmath}
\usepackage{amssymb}
\usepackage{graphics,epsfig}
\newtheorem{theorem}{Theorem}[section]
\newtheorem{lemma}[theorem]{Lemma}
\newtheorem{proposition}[theorem]{Proposition}
\newtheorem{corollary}[theorem]{Corollary}

\newtheorem{definition}[theorem]{Definition}
\newtheorem{example}[theorem]{Example}

\newtheorem{remark}[theorem]{Remark}

\newcommand{\PGLn}[1]{\mbox{PGL} (#1, \mathbb R)}
\newcommand{\aff}[1]{\mbox{Aff} (#1,\mathbb R)}
\newcommand{\rp}[1]{\mathbb {R P}^{#1}}
\newcommand{\subinf}[1]{\inf_{#1 \in G}}

\newcommand{\lemref}[1]{Lemma~\ref{#1}}
\newcommand{\propref}[1]{Proposition~\ref{#1}}

\begin{document}
\title{Quasi-homogeneous domains and convex affine manifolds}
\author{Kyeonghee Jo}
\begin{abstract}
In this article, we study convex affine domains which can cover a
compact affine manifold.

For this purpose, we first show that every strictly convex
quasi-homogeneous projective domain has at least $C^1$ boundary
and it is an ellipsoid if its boundary is twice differentiable.
And then we show that an n-dimensional paraboloid is the only
strictly convex quasi-homogeneous affine domain in $\mathbb R^n$
up to affine equivalence. Furthermore we prove that if a strictly
convex quasi-homogeneous projective domain is $C^{\alpha}$ on an
open subset of its boundary, then it is $C^{\alpha}$ everywhere.

Using this fact and the properties of asymptotic cones we find all
possible shapes for developing images of compact convex affine
manifolds with dimension $\leq 4$.
\end{abstract}
\subjclass[2000]{52A20, 57M50} \keywords{quasi-homogeneous,
divisible, strictly convex, projective manifold, affine manifold}
\address{Applied Mathematics, Sejong University,
143-747, Seoul, Korea} \email{khjo@sejong.ac.kr}
\thanks {This work was supported by Korea Research Foundation Grant (KRF-2002-070-C00010).}
\maketitle

\section{Introduction}
A convex affine n-manifold is a quotient of a convex affine domain
$\Omega$ in $\mathbb R ^n$ by a discrete subgroup $\Gamma$ of
$\aff{n}$ acting on $\Omega$ properly discontinuously and freely.
We say that a convex affine manifold is complete if $\Omega =
\mathbb R ^n$ and radiant if $\Omega$ is a cone. Particularly if
$M$ is a compact convex affine manifold then $\Omega$ becomes a
divisible affine domain, i.e., its automorphism group contains a
cocompact discrete subgroup acting properly. So the study of
compact convex affine manifolds is equivalent to that of
divisible convex affine domains and their automorphism groups.

Compact convex affine 3-manifolds were classified completely: The
complete case was treated by D. Fried and W. Goldman in \cite{FG}
and the radiant case was treated by T. Barbot in \cite{Bt}. The
remaining cases were done by H. Cho in \cite{Cho}, where she
classified them by finding all 3-dimensional divisible convex
affine domains which are not cones independently.

Recently, Y. Benoist proved in \cite{Bn2} that for any properly
convex cone $\Omega \subset \mathbb R^n$ and a discrete subgroup
$\Gamma$ of Aff$(n,\mathbb R)$ which divides $\Omega$, one of the
following three cases holds true:
\begin{enumerate}
\item[\rm(i)] $\Omega$ is a product,
\item[\rm(ii)] $\Omega$ is a symmetric cone,
\item[\rm(iii)] $\Gamma$ is Zariski dense in $\mbox{GL}(n, \mathbb R)$.
\end{enumerate}

One of the oldest results in this subject is that any divisible
convex affine domain $\Omega$ must be a cone if it does not
contain any complete line, which was proved by Vey in \cite{V3}.
If $\Omega$ has a complete line it is isomorphic to $\mathbb R^k
\times \Omega'$ for some $(n-k)$-dimensional convex affine domain
$\Omega'$ which does not contain any complete line. We see that
the action of $\Gamma$ on $\Omega$ induces a quotient affine
action $\tilde{\Gamma}$ on $\Omega'$. In general, $\tilde{\Gamma}$
may not be discrete any more. But $\Omega'$ still has a compact
quotient by the action of $\tilde{\Gamma}$ and thus becomes just a
quasi-homogeneous convex affine domain, i.e., there is a compact
subset $K$ of $\Omega'$ such that $\tilde{\Gamma}K = \Omega'$. So
any divisible convex affine domain is isomorphic to the product of
$\mathbb R^k$ and a quasi-homogeneous convex affine domain which
does not contain any complete line. A good example is $\mathbb R
\times \{(x,y)\in \mathbb R^2 \,|\,y>x^2\}$ (see \cite{G4} for
details).

For this reason we study quasi-homogeneous convex affine domains
to understand compact convex affine manifolds. In this article we
first show that an $n$-dimensional paraboloid $\{(x_1, x_2, \dots,
x_n)\in \mathbb R^n\,|\,x_n>x_1^2+x_2^2+ \dots+x_{n-1}^2\}$
 is the only strictly convex quasi-homogeneous affine
domain in $\mathbb R^n$ up to affine equivalence, in contrast with
the fact that there are infinitely many quasi-homogeneous strictly
convex projective domains (see \cite{G1,Ko,VK}). To prove this we
investigate strictly convex quasi-homogeneous projective domains,
 since strictly convex quasi-homogeneous affine domains are a special class of them.
In general, any affine domain of  $\mathbb R^n$ can be
 considered as a projective domain via the well-known
 equivariant embedding from $(\mathbb R^n, \aff{n}) $ into   $(\mathbb{RP}^n,  \mbox{PGL}(n+1,\mathbb R))$.
  So we often look at affine domains in the projective space. This enables us to use Benz\'{e}cri's technique developed
in convex projective domain theory \cite{B}.

Actually in section \ref{SC} we show that if $\Omega$ is a
strictly convex quasi-homogeneous projective domain, then
\begin{enumerate}
\item[\rm (i)] $\partial \Omega$ is at least
$C^1$,
\item[\rm (ii)] $\Omega$ is an ellipsoid if $\partial \Omega$ is
twice differentiabe,
\item[\rm (iii)] if $\partial \Omega$ is $C^{\alpha}$ on an open subset of $\partial \Omega$, then
$\partial \Omega$ is $C^{\alpha}$ everywhere.
\end{enumerate}
 This implies that either any strictly convex projective domain
which covers a compact projectively flat manifold is an ellipsoid
or the boundary fails to be twice differentiable on a dense
subset, which is a generalization of the 2-dimensional result of
Kuiper \cite{Ku}.

For divisible case, a similar result has been proved independently
by Y. Benoist in \cite{Bn1} and he claimed in \cite{Bn2} that any
quasi-homogeneous strictly convex projective domain which is not
an ellipsoid, has a discrete automorphism group. We will also show
the claim in section \ref{strict3}.

Another earlier result related to this subject is that if $\Omega$
is a strictly convex (in the sense that Hessian is positive
definite) domain with $C^3$ boundary which has a cofinite volume
discrete subgroup action, then $\Omega$ is an ellipsoid \cite{CV}.

To conclude that any strictly convex quasi-homogeneous affine
domain is affinely equivalent to a paraboloid, we show in section
\ref{strict2} that its boundary is homogeneous and thus
$C^{\infty}$. (Then the domain is projectively equivalent to an
ellipsoid by the above result (ii) about quasi-homogeneous
strictly convex projective domain and this implies that it is
affinely equivalent to a paraboloid.)
 To do this we first prove that every strictly convex domain
has an invariant direction under the action of its affine
automorphism group, that is, it contains a half line which is
invariant under the action of the linear part of its automorphism
group. Generally, we will show in section \ref{AC} that every
convex quasi-homogeneous affine domain $\Omega$ is foliated by
cosets of the {\it asymptotic cone} which is the maximal cone
contained in $\Omega$ and invariant under the action of the linear
part of the its automorphism group. Using this asymptotic cone we
find all the possible quasi-homogeneous convex affine domains with
dimension $\leq 3$. Then we find all the possible shapes for a
universal covering space of compact convex affine manifold with
dimension $\leq 4$. Finally, we give a partial result to the
Markus conjecture which says that ``A compact affine manifold is
complete if and only if it has parallel volume.'' This conjecture
was proved under some additional conditions about holonomy group
by J. Smillie, D. Fried, W. Goldman, M. Hirsch, Carri\`{e}re and
so on, see \cite{C1}, \cite{C2}, \cite{F} \cite{FGH}, \cite{GH1},
\cite{GH}, \cite{Ma}, \cite{Sm}. In this article we show that any
compact convex affine manifold with dimension $\leq 4$ is complete
if it has parallel volume.

\section{Preliminaries}\label{pre}
We present here some of the basic materials from Benz\'{e}cri's
convex domain theory that we will need later. We begin with some
definitions.

The real projective space $\rp{n}$ is the quotient space of
$\mathbb R^{n+1} \setminus \{0\}$ by the action of $\mathbb
R^*=\mathbb R \setminus \{0\}$. In an affine space, we usually
denote the affine subspace
 generated by a subset $A$ by $\langle A\rangle$. So we will use the same
 notation for a subset of $\rp{n}$, i.e., for each subset $B$ of
 $\rp{n}$ $\langle B \rangle$ means the projectivization of the affine subspace
 generated by $\pi^{-1}(B)$ in  $\mathbb R^{n+1}$, where $\pi$
 is the quotient map from  $\mathbb R^{n+1} \setminus \{0\}$ onto $\rp{n}$ and
 we will call $\langle B \rangle$ the $support$ of $B$.

A $quasi$-$homogeneous$ affine (respectively, projective) domain is an an
open subset $\Omega$ of $\mathbb R^n$ (respectively, $\rp{n} $) which has
a compact subset $K \subset \Omega$
 and a subgroup $G$ of Aut$(\Omega)$ such that $GK=\Omega$, where
 Aut$(\Omega)$ is a subgroup of $\mbox{Aff}(n,\mathbb R)$ (respectively, $\PGLn{n+1}$)
 consisting of all affine (respectively, projective) transformations preserving
 $\Omega$. Sometimes we say that $G$ acts on $\Omega$ $syndetically$.
 It follows that both
 homogeneous domains and divisible domains are quasi-homogeneous.

 As stated in the introduction, we can consider an affine
domain as a projective domain and it is obvious that every
quasi-homogeneous affine domain is a quasi-homogeneous projective
domain. But a quasi-homogeneous projective domain is not a
quasi-homogeneous affine domain even if it is contained in an
affine patch of a projective space. For examples, any
quasi-homogeneous strictly convex projective domain which is not
an ellipsoid, is not a quasi-homogeneous affine domain (this will
be proved in section 7) and in fact there exist infinitely many
such quasi-homogeneous strictly convex projective domains. (See
\cite{G1} and \cite{VK} for the 2-dimensional case and \cite{Bn0}
for arbitrary dimensional cases.)

For an affine domain $\Omega \subset \mathbb R^n$, we denote the
group of all affine transformations preserving $\Omega$ by
$\mbox{Aut}_{\text{aff}}(\Omega)$ to distinguish it from the group
of all projective transformations preserving $\Omega$. Note that
$\mbox{Aut}_{\text{aff}}(\Omega)$ is a closed subgroup of
Aut$(\Omega)$. Also usually we denote the boundary of a domain
$\Omega$ by $\partial \Omega$, but sometimes we will use the
$\partial_a \Omega$ and $\partial_p \Omega$ when it is necessary
to avoid ambiguity and call them an affine boundary and a
projective boundary of $\Omega$, respectively.
 Note that $\partial_a \Omega$ is a subset of $\partial_p \Omega$ and
in fact $\partial_a \Omega=\mathbb R^n \cap \partial_p \Omega$.

\begin{figure}[h!]
\begin{center}
 \epsfig{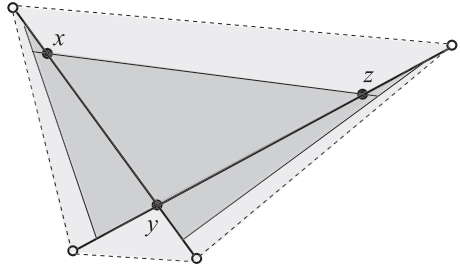}
 \caption{Transitivity of $\sim$
\label{fig:angle}}
\end{center}
\end{figure}

For a convex projective domain $\Omega \subset \rp{n}$, we define
a relation $ \sim$ on $\overline \Omega$ by defining $x \sim y$ if
either $x = y$ or $x \neq y$ and  $\overline \Omega$ has an open
line segment $l$ containing both $x$ and $y$. It follows from the
convexity of $\Omega$ that $\sim$ is an equivalence relation : its
symmetry is obvious by definition. See Fig. \ref{fig:angle} for
its transitivity.

 An equivalence class with respect to $\sim$ is called a $face$ of
 $\Omega$ and we will call the closure of a face a $closed$ $face$ of $\Omega$.
  Note that a face is relatively open in its support and $\overline\Omega$ is a disjoint union of all faces.
  A convex domain $\Omega$ in $\rp{n}$ is called $properly$ $convex$ if there is no
non-constant projective map of $\mathbb R$ into $\Omega$ and
$strictly$ $convex$ if $\partial \Omega$ contains no line segment.
From these definitions we see that any strictly convex domain is a
properly convex domain.

It is clear that the intersection of any family of convex sets is
again convex. Therefore, for any subset $A$ there is a smallest
convex set containing $A$, namely, the intersection of all convex
sets containing $A$. This convex set is called the convex set
spanned by $A$, or the convex hull of $A$, and is denoted by
$CH(A)$.

An affine subspace $\mathcal S$ is called a \emph{ supporting
subspace} of a convex domain $\Omega$ in $\rp{n}$ if following
conditions are satisfied :
\begin{enumerate}
\item [\rm (i)]$\mathcal S\cap \overline\Omega \neq \emptyset $,
and
\item [\rm (ii)]for any face $F$ of $ \Omega $ such $F\cap \mathcal S\neq
\emptyset$, $F$ is entirely contained in $\mathcal S$.
\end{enumerate}
We see immediately that the support of a face is always a
supporting subspace and any hyperplane $H$ tangent to $\partial
\Omega$ is a supporting subspace of $\Omega$.
\begin{example}
\item [\rm (i)] Let $\Omega=\{(x,y,z) \in \mathbb
R^3\,|\, y>x^2, z>0\}$.
 Then the $z$-axis is the support of a face $\{(0,0,z) \in \mathbb
R^3\,|\,  z>0\}$,
 the $y$-axis is neither a support of any face nor a supporting subspace of $\Omega$, and the $x$-axis is not a
 support of any face but a supporting subspace of $\Omega$.

\item [\rm (ii)] Let $\Omega$ be the standard simplex in $\mathbb R^{n}$, that is, its vertices are
 $\{0,e_1,\dots, e_n \}$. Then $\{(x_1,\dots,x_n)\in \mathbb R^{n} \,|\, x_n=1\}$ is a supporting subspace of $\Omega$.
 (Note that this is not a quasi-homogeneous affine domain but a quasi-homogeneous projective domain.)
\end{example}

The following definitions are equivalent to Benz\'{e}cri's
original definitions in \cite{B}.

\begin{definition}\label{def-Ben}
\mbox{}
\begin{enumerate}
\item [\rm (i)] \emph{ Let $\Omega$ be a properly convex domain in $\rp{n}$ and let $\Omega_1$ and
$\Omega_2$ be convex domains in $\langle \Omega_1 \rangle$ and
$\langle \Omega_2 \rangle$ respectively. $\Omega$ is called a}
convex sum \emph{of $\Omega_1$ and $\Omega_2$, which will be
denoted by $\Omega = \Omega_1 + \Omega_2$, if $\langle \Omega_1
\rangle \cap \langle \Omega_2 \rangle =\emptyset$ and $\Omega$ is
the union of all open line segments joining points in $\Omega_1$
to points in $\Omega_2$.}
\item [\rm (ii)]\emph{A $k$-dimensional face $F$ of an $n$-dimensional convex domain $\Omega$ is
called} conic \emph{if there exist $n-k$ supporting hyperplanes
$H_1, H_2, \dots,H_{n-k}$ such that}
\begin{equation*}
H_1 \gneq H_1 \cap H_2 \gneq \dots \gneq H_1 \cap \dots \cap
H_{n-k} = \langle F \rangle.
\end{equation*}
\item [\rm (iii)] \emph{Let $\Omega$ be a properly convex domain in $\rp{n}$.
We say that $\Omega$ has an $\emph{osculating}$}  ellipsoid at
\emph{$p \in
\partial \Omega$ if there exist a suitable affine chart and a basis
such that the local boundary equation on some neighborhood of
$p=(0,\dots,0)$ is expressed by $x_n = f(x_1,\dots, x_{n-1})$ and}
\begin{equation*} \lim_{(x_{1},\dots,x_{n-1})\to 0}
\frac{f(x_{1},\dots,x_{n-1})}{{x_1}^2+ \dots +{x_{n-1}}^2}=1.
\end{equation*}
\end{enumerate}
\end{definition}
\begin{remark}
From Definition \ref{def-Ben} (ii), any (n-1)-dimensional face is
conic.
\end{remark}
\begin{example}
\item [\rm (i)] Let $\Omega$ be a simplex. Then every face of $\Omega$ is conic.
\item [\rm (ii)] Let $\Omega=\{(x,y,z) \in \mathbb
R^3\,|\, y>x^2, z>0\}$.
 Then it is a quasi-homogeneous affine (and so projective) domain and the subset
  $\{(0,0,z) \in \mathbb R^3\,|\,  z>0\}$ is a non-conic face of $\Omega$.
\end{example}
\subsection{Spaces of convex bodies in $\rp{n}$}
In 1960, Benz\'{e}cri developed projective convex body theory in
\cite{B}. We recall here some important results among them, which
can be also found in \cite {G3}. We will follow Goldman's notation
\cite {G3}.

First, we call a subset $K \subset \rp{n}$ a $convex$ $body$ if
$K$ is the closure of a properly convex domain of $\rp{n}$.
  Equivalently, $K$ is a convex body if $K$ is a closed convex
subset with nonempty interior of $\rp{n}$ which does not contain
any complete line.

Let $\mathcal{C}(n)$ denote the set of all convex bodies in
$\rp{n}$, with the topology induced from the Hausdorff metric on
the set of all closed subsets of $\rp{n}$. Let
$$\mathcal{C_*}(n)=\{ (K,x) \in  \mathcal{C}(n)\times \rp{n} \,|\,
x \in \mbox{int}(K) \}$$ be the corresponding set of pointed
convex bodies with topology induced from the product topology on $
\mathcal{C}(n)\times \rp{n}$.

\begin{theorem}[Benz\'{e}cri] $\PGLn{n+1}$ acts properly and
syndetically on $\mathcal{C_*}(n)$. In particular the quotient
$\mathcal{C_*}(n)/\PGLn{n+1}$ is a compact Hausdorff space.
\end{theorem}
\begin{remark} While the the quotient
$\mathcal{C_*}(n)/\PGLn{n+1}$ is Hausdorff, the space of
equivalence classes of convex bodies $\mathcal{C}(n)/\PGLn{n+1}$
is not Hausdorff. Some basic examples can be found in \cite{B} and
in C.16 of \cite {G3}.
\end{remark}

Using  the above theorem, he proved the following important result.
\begin{theorem}[Benz\'{e}cri]\label{thm-benz}
Let $\Omega \subset \rp{n}$ is a quasi-homogeneous properly convex
domain. Then the corresponding point $[\overline{\Omega}] \in
\mathcal{C}(n)/\PGLn{n+1}$ is closed, that is, if $g_{n}\Omega $
converges to some properly convex domain $\Omega' \subset \rp{n}$
for some $\{g_{n}\}\subset \PGLn{n+1}$, then $\Omega $ is
projectively equivalent to $\Omega'$.
\end{theorem}

\begin{theorem}[Benz\'{e}cri]\label{prop-benz}
\begin{enumerate}
\item [\rm(i)]
Let $\Omega$ be a properly convex projective domain in $\rp{n}$
and $F$ be a conic face of $\Omega$. Then there exist a projective subspace $L$ of $\rp{n}$ and
projective automorphisms $\{h_i\}$ of $\rp{n}$
such that $\{h_i\Omega \}$ converges to $F+B$ for some properly
convex domain $B$ in $L$.
\item[\rm(ii)]
If a properly convex domain $\Omega$ has an osculating ellipsoid
$Q$, then there exists a sequence $\{g_{n}\}\subset \PGLn{n+1}$
such that $g_{n}\Omega $ converges to $Q$.
\item[\rm(iii)]
Let $\Omega = {\Omega}_1 + {\Omega}_2$. Then $\Omega$ is
quasi-homogeneous(respectively, homogeneous) if and only if ${\Omega}_i$
is quasi-homogeneous(respectively, homogeneous) for each i.
\item[\rm(iv)]
 Let $\Omega$ be a quasi-homogeneous properly convex domain in $\rp{n}$ and $L$ a
linear subspace of $\rp{n}$ of dimension $r$ such that $L \cap
\Omega$ has a conic face $F$. Then there exists a section which is
projectively equivalent to an $r$-dimensional properly convex
domain $F+B$ for some suitable properly convex domain $B$ of
dimension $r- (\dim(F)+1)$.

\end{enumerate}
Here a section means an intersection of a projective subspace and $\Omega$.
\end{theorem}

\begin{corollary}\label{cor-con}
Any quasi-homogeneous properly convex projective domain which is
not strictly convex has a triangle as a 2-dimensional section.
\end{corollary}
{\bf Proof. } Let $\Omega$ be a domain satisfying the hypothesis.
Then $\partial \Omega$ contains a line segment $l$. Considering a
2-dimensional plane $P$ which intersects $\Omega$ and contains
$l$, we see that $l$ is a conic face of $P\cap \Omega$. By Theorem
\ref{prop-benz} (iv), we get the desired result.
 $\qquad \square $

\section{Aut($\Omega$) and its limit projective transformations}
  Any properly convex projective domain
is projectively equivalent to a bounded convex domain in an affine
space. So we can define on any properly convex projective domain a
complete continuous metric which is invariant under the action of
Aut$(\Omega)$. This metric is called the Hilbert metric and
defined as follows.
\begin{definition}
\emph{Let $\Omega$ be a bounded convex domain in $\mathbb R^n$.
For any two different points $p_1, p_2 \in \Omega$, we define
$\mbox{d}_H(p_1, p_2)$ to be the logarithm of the absolute value
of the cross ratio of $(s_1, s_2, p_1, p_2)$, where $s_1$ and
$s_2$ are the points in which the line $\overline {p_1p_2}$
intersects $\partial \Omega$. For $p_1 = p_2$, we define
$\mbox{d}_H(p_1, p_2) =0$.}
\end{definition}

Since $\mbox{PM}(n+1,\mathbb R)$, which is the projectivization of
the group of all (n+1) by (n+1) matrices, is a compactification of
$\PGLn{n+1}$, any infinite sequence of non singular projective
transformations contains a convergent subsequence. Note that the
limit projective transformation may be singular. For a singular
projective transformation $g$ we will denote the projectivization
of the kernel and range of $g$ by $K(g)$ and $R(g)$. Then $g$ maps
$\rp{n} \setminus K(g)$ onto $R(g)$ and the images of any compact
set in $\rp{n} \setminus K(g)$ under the convergent sequence
${g_i}$, converge uniformly to the images under the limit
transformation $g$ of ${g_i}$ (see \cite {B} and \cite{Hyuk}).

\begin{lemma}\label{lem-Ben}Let $\Omega$ be a properly convex domain
in $\rp{n}$ and let $\{g_i \} \subset \mbox{Aut} ( \Omega ) $ be a
sequence converging to a singular projective transformation $g \in
\mbox{PM} (n+1,\mathbb R)$. Then,
\begin{enumerate}
  \item [\rm(i)]$K(g)$ is a supporting subspace of $\Omega$,
  \item [\rm(ii)]$R(g)=\langle R(g) \cap \partial \Omega \rangle$ and
  $R(g)$ is a supporting subspace of $\Omega$.
\end{enumerate}
\end{lemma}
{\bf Proof. }In \cite[p310--311]{B}, Benz\'{e}cri proved that $K(
g)$ is a supporting subspace and $R(g)=\langle R(g)\cap \partial
\Omega \rangle$. So it suffices to show that $R(g)$ is a
supporting subspace of $\Omega$.
 Suppose not. We
 can choose two points $\{x,y\}$ which are contained in $\partial
 \Omega \cap R(g)^c$ such that $\overline{xy} \cap R(g)$ is a
 point, that is, $\overline{xy}$ meets $R(g)$ transversally. Then we
 can find a small neighborhood $U$ of $R(g) \cap \partial
 \Omega$ in $\rp{n}$ which does not contain any $\varepsilon$-ball
 with respect to the Hilbert metric on $\Omega$. But any compact
 $\varepsilon$-ball neighborhood $B(z,\varepsilon)$ of a point $z$
 in $\Omega$ is mapped into $U \cap \overline \Omega$ under the $g_i$ action
 for sufficiently large $i$, since $g_i(B(z,\varepsilon))$
 converges uniformly to $g(B(z,\varepsilon)) \subset R(g)\cap \partial
 \Omega$. This contradiction completes the proof.
 $\qquad \square $
\begin{remark}
\mbox{}

\begin{enumerate}
\item
Note that neither $K(g)$ nor $R(g)$ intersects $\Omega$ since
$\emph{dim}K(g)<n$ and $\emph{dim}R(g)<n$ in this situation.
\item It follows from (ii) of the above lemma that $R(g)\cap \partial \Omega$ is a closed face of $\Omega$. In fact, we
can show that any supporting subspace $L$ satisfying $L=\langle L \cap \partial \Omega \rangle$
is always a support of a face.
\end{enumerate}
\end{remark}

\begin{lemma} \label{lem-saillant}
Let $\Omega$ be a quasi-homogeneous properly convex domain in
$\rp{n}$and $G$ a subgroup of $\mbox{Aut} ( \Omega ) $ acting on
$\Omega$ syndetically. Then we get the following :
\begin{enumerate}
\item [\rm(i)] For each point $p\in \partial \Omega$,
there exists a sequence $\{g_i\}\subset G$ and $x\in \Omega$ such
that $g_{i}(x)$ converges to $p$.
\item [\rm(ii)]For any accumulation point $g$ of $\{g_{i}\}$, $R(g)$ is the
projective subspace generated by the face containing $p$.
\end{enumerate}
\end{lemma}

{\bf Proof. }
\begin{enumerate}
\item [\rm(i)]
 Choose a
 straight line segment $\ell$ in $\overline {\Omega}$ with initial
 point in $\Omega$ and endpoint $p$. Let $K$ be a compact generating
 domain of $\Omega$. Then $\ell$ intersects an infinite series of
 different images of $K$ ; $g_{1}(K),g_{2}(K),\dots .$ Choose $x_{i}\in \ell\cap
 g_{i}(K)$ for all $i$. Then we have $g_{i}^{-1}(x_{i})\in K$ and thus we
 can choose $x\in K$ which is one of the accumulation points of
 $\{g_{i}^{-1}(x_{i})\}$ since $K$ is compact.  We may
 assume $\{g_{i}\}$ converges to $g$. Since
 $d_H (g_{i}(x),x_{i})=d_H (x,g_{i}^{-1}(x_{i}))$ converges to $0$,
 $d(g_{i}(x),x_{i})$ also converges to $0$, where $d_H$ is a
 Hilbert metric and $d$ is the standard metric of $\rp{n}$.
 We see that $d(g_{i}(x),p)\leq d(g_{i}(x),x_{i})+ d(x_{i},p)$ and this
implies that $g_{i}(x)$ converges to $p$.
\item [\rm(ii)] Obviously $g$ is a singular projective transformation.
Let $F$ be the face of
$\Omega$ which contains $p$.
Then by \lemref{lem-Ben} $R(g)$ contains $F$ and $R(g)=\langle R(g) \cap \partial \Omega \rangle$.
Since $p$ must be an interior point of $R(g) \cap \partial \Omega$, we get
 $R(g)=\langle F \rangle$.

\end{enumerate}
 $\qquad \square $

\begin{lemma}\label{lem-kernel}
Let $\{f_{i}\}$ be a sequence in $\aff{n} (<\mbox{PGL}(n+1,\mathbb
R))$. Suppose that $f_i$ converges to $f\in  \mbox{PM}(n+1,\mathbb
R)$ with $R(f)\cap \mathbb R^n \neq \emptyset$. Then $K(f)\cap
\mathbb R^n=\emptyset$.
\end{lemma}
{\bf Proof. }Suppose $K(f) \cap \mathbb R^n \neq \emptyset$ and
choose $x\in K(f)\cap \mathbb R^n$.  For any one-dimensional
affine subspace
 $L$ such that $L\cap K(f)=\{x\}$,
$L\setminus \{x\}$ is mapped into a point in $R(f)$ by $f$. Since
every one-dimensional affine subspace through $x$ meets $\partial
\mathbb R^n =\rp{n-1}_{\infty}$, for any $p\in R(f)$ there exist
$q\in \partial \mathbb R^n \backslash K(f)$ such that $f_{i}(q)$
converges to $p=f(q)$. Now we choose any point $y$ in $\mathbb R^n
\setminus K(f)$ such that $f(y)\in R(f)\cap \mathbb R^n$. Then the
one-dimensional affine subspace
 $l_{xy}$ containing $x$ and $y$ intersects
$\partial \mathbb R^n$ at a point $z$. Since $l_{xy}\cap
K(f)=\{x\}$, $z$ is not contained in $K(g)$ and so $f_i(z)$
converges to $f(z)=f(y)$.
 But
this is impossible since $f_i(z) \in f_{i}(\partial \mathbb R^n
)=\partial \mathbb R^n$ for all $i$ and $f(y) \in \mathbb R^n$.
This completes the proof.
 $\qquad \square $

\begin{lemma}\label{lem-bounded} Let $\Omega$ be a quasi-homogeneous properly convex affine domain in $\mathbb R^n$.
Then $\partial_a \Omega$ does not have any bounded face with
dimension $k>0$.
\end{lemma}
{\bf Proof. }Let $L$ be a supporting subspace of $\dim k$ such
that $\mathcal S=\partial_a \Omega \cap L$ is a bounded closed
subset of $L$ with nonempty interior. By Lemma
~\ref{lem-saillant}, there exists a sequence $\{f_i\}\subset
\mbox{Aut}_{\text{aff}} (\Omega)$ which converges to $f\in
\mbox{PM}(n+1,\mathbb R)$ with $R(f)=L$ since int$\mathcal S$ is a
face of $\Omega$. Then Lemma ~\ref{lem-kernel} implies that
$K(f)\subset
\partial \mathbb R^n$ and $K(f)\neq \partial \mathbb R^n$ since $k>0$.
Note that $f_i(\partial \mathbb R^n)=\partial \mathbb R^n$ for
each $i$ and this implies $f(\partial \mathbb R^n \setminus
K(f))\subset \partial \mathbb R^n$.

Let $K$ be the linear subspace of all points with direction in
$K(f)$. Then $K$ is $(n-k)$-dimensional. Choose any
$k$-dimensional linear subspace $N$ which is transversal to $K$
and consider a projection $p$ from $\mathbb R^n$ to $N$ satisfying
$p(x+k)=x$ for any $x \in \mathbb R^n$ and any $k \in K$. Suppose
$p(\Omega)$ is not bounded in $N$. Then there exist sequences
$\{x_n\} \subset N$ and $\{y_n\} \subset \Omega$ such that $y_n
\in x_n+K$ and $x_n$ converges to a point $x_{\infty}$ in
$\overline N \cap
\partial \mathbb R^n \subset \partial \mathbb R^n \setminus K(f)$,
where $\overline N$ means the closure of $N$ in $\rp{n}$.  Observe
that $f(x_n) = f(y_n)$ since $y_n$ lies in the projection of
$\pi^{-1}(\langle x_n,K(f) \rangle)$. But $f(x_n)$ converges to a
point $f(x_{\infty})$ in $\overline L \cap \partial \mathbb R^n $
and $f(y_n)$ is a point $\mathcal S=L \cap \partial \Omega$. This
is a contradiction. So $p(\Omega)$ must be bounded in $N$. This
implies that $\Omega$ must lie between two parallel hyperplanes.
But it is known that any quasi-homogeneous affine domain can not
lie between parallel hyperplanes (see \cite {Park}).
 $\qquad \square $

For a properly convex open cone Vey showed in \cite{V3} that
$\mbox{Aut}_{\text{aff}}(\Omega)$ acts properly on $\Omega$ using
the properties of the characteristic function of $\Omega$. In fact
we can prove the following.

 \begin{proposition}\label{prop-acts}
      Let $\Omega $ be a properly convex domain in
     $\rp{n}$ and $G$ be a closed subgroup of $\mbox{Aut} (\Omega )$. Then $G$ acts properly on $\Omega$.
     \end{proposition}
     {\bf Proof. } Let $K$ and $L$ be compact subsets of $\Omega$.
     Define $C(K, L) \subset \mbox{Aut} (\Omega)$ by
     \[C(K, L)=\{g\in G \,|\, gK\cap L\neq \emptyset\}.\]
     Let $\{g_n\}\subset C(K, L)$ be a sequence converging to $g\in \mbox{PM}(n+1,\mathbb  R)$.
     By the definition of $C(K, L)$, there exists $x_n \in K$ such that $g_n
x_n \in L$ for each $n$. Since $K$ and $L$ are compact, we may
assume that $x_n$ converges to $x\in K$ and $g_n x_n $ converges
to $y\in L$ by choosing a subsequence if necessary. Then,
 \begin{align*}
\mbox{d}_{H}(g_{n}x,y)&\leq
\mbox{d}_{H}(g_{n}x,g_{n}x_{n})+\mbox{d}_{H}(g_{n}x_{n},y)\\
&=\mbox{d}_{H}(x,x_{n})+\mbox{d}_{H}(g_{n}x_{n},y), \end{align*}
 \begin{equation*}
\lim_{n\to \infty}\mbox{d}_{H}(g_{n}x,y)=0,\mbox{ i.e. }\lim_{n\to
\infty}g_{n}x=y.
\end{equation*}

Suppose $g$ is singular. By \lemref{lem-Ben} (i), $x\notin K(g)$
and thus $gx=y$ is contained in $R(g)$. This contradicts
\lemref{lem-Ben} (ii). Therefore $g$ is nonsingular and so $g\in
C(K, L)$. This implies that $\overline {C(K, L)} =C(K, L) $ in
$\mbox{PM}(n+1,\mathbb R)$ and we conclude that $C(K, L)$ is
compact.
 $\qquad \square $

The following lemma is actually an immediate corollary of the above proposition,
but we give a direct proof here.
\begin{lemma}\label{closed} Let $\Omega$ be a properly convex
projective domain in $\mathbb {RP}^n$ and $G$ be a closed subgroup
of $\mbox{Aut}(\Omega)<\mbox{PGL}(n+1,\mathbb R)$. Then each
$G$-orbit is a closed subset of $\Omega$.
\end{lemma}
{\bf Proof. } Let $x \in \Omega$ and $y\in \overline {Gx} \cap
\Omega$. Choose a sequence $\{g_n\}$ of automorphisms such that
$\lim_{n\to \infty}g_n x=y$. Then there exists a subsequence
$\{g_{n_k}\}$ of $\{g_n\}$ converging to some $g\in
\mbox{PM}(n+1,\mathbb R)$ by compactness of $\mbox{PM}(n+1,\mathbb
R)$. Suppose $g$ is singular. Then by Lemma \ref{lem-Ben} we know
$x\notin K(g)$ and this implies $y\in R(g)$. This contradicts
Lemma \ref{lem-Ben}. So $g$ can not be singular, that is, $g\in
\mbox{Aut}(\Omega)$. Thus $y=gx\in Gx$ since $G$ is closed in
Aut$(\Omega)$.
 $\qquad \square $

The following theorem is well-known and a proof is given in
Theorem 6.5.1. of \cite{R2}.
\begin{theorem}\label{metric}Let $\Gamma$ be a group of isometries of
a metric space $X$ and $d_{\Gamma}$ the orbit space distance
function
\begin{equation*}
\emph{d}_{\Gamma}:X/\Gamma \times X/\Gamma \to \mathbb R
\end{equation*}
defined by the formula $$\emph{d}_{\Gamma}(\Gamma x,\Gamma
y)=\emph{dist}(\Gamma x,\Gamma y)= \inf\{d(z_1,z_2)\,|\,z_1 \in
\Gamma x \mbox{ and } z_2 \in \Gamma y \}.$$ Then
$\mbox{d}_{\Gamma}$ is a metric on $X/\Gamma $ if and only if each
$\Gamma$-orbit is a closed subset of $X$.

\end{theorem}
 By Lemma \ref{closed} and Theorem \ref{metric}, we obtain a metric $\mbox{d}_{G}$ on
the orbit space $G\setminus \Omega $ using $\mbox{d}_{H}$ as
follows:
\begin{equation}\label{eq-metric}
 \mbox{d}_{G}(Gx,Gy)=\subinf{g,g'} \mbox{d}_{H}(gx,g'y).
\end{equation}
Since $\mbox{d}_{H}$ is invariant by projective automorphism, we
have
\begin{equation*}
\mbox{d}_{G}(Gx,Gy)=\subinf{g} \mbox{d}_{H}(x,gy).
\end{equation*}
Up to now we have shown that the Hilbert metric on a properly
convex domain $\Omega$ descends to the orbit space $G\setminus
\Omega $ if $G$ is a closed subgroup of Aut$(\Omega)$. In
addition, we can show that the distance between the two orbits
$Gx$ and $Gy$ is realized by the Hilbert distance between two
points $x_0 \in Gx$ and $y_0 \in Gy$ in $\Omega$.
\begin{lemma}\label{lem-exist} Let $G$ be a closed subgroup of $\PGLn{n+1}$. Then for
each pair $(Gx,Gy)\in G\setminus \Omega $, there exists $g_{x,y}
\in G$ such that $\mbox{d}_{G}(Gx,Gy)=\mbox{d}_{H}(x,g_{x,y}y)$.
\end{lemma}
{\bf Proof. } By \eqref{eq-metric}, there exists a sequence $\{g_n
\} \subset G$ such that $\mbox{d}_{H}(x,g_{n}y)$ converges to
$\mbox{d}_{G}(Gx,Gy)$. For any $r > \mbox{d}_{G}(Gx,Gy)$, there
exists $N>0$ such that $g_{n}y \in \mbox{B}(x,r)$ whenever $n\geq
N$. We may assume that $\{g_{n}\}$ converge to a projective
transformation $g\in \mbox{PM}(n+1,\mathbb R)$ by taking a
subsequence if necessary. If $g$ is singular then $g_{n}y$
converges to some point in $\partial \Omega$. This is a
contradiction. Therefore $\{g_{n}\}$ converges to $g_{x,y} \in G$
and so we get
\begin{equation*}
\mbox{d}_{H}(x,g_{x,y}y)=\lim_{n\to \infty} \mbox{d}_{H}(x,g_{n}y)
\end{equation*}
since $g_{n}y$ converges to $g_{0}y\in \Omega $.
 $\qquad \square $


\section{Asymptotic cone}\label{AC}

As a special subclass of quasi-homogeneous convex projective
domains, quasi-homogeneous convex affine domains have some
important properties which distinguish them from other projective
domains. They contain a cone invariant under the action of linear
parts of their automorphism groups, which is called an asymptotic
cone. This terminology was originally introduced by Vey in
\cite{V3}.

\begin{definition}

 \emph{Let $\Omega$ be a convex domain in $\mathbb
R^n$. The} asymptotic cone \emph{of $\Omega$ is defined as
follows:}
\begin{equation*}
\emph{AC} (\Omega)=\{u\in \mathbb R^n\,|\, x+tu\in \Omega, \mbox{
for all } x\in \Omega, t\geq 0\}.
\end{equation*}
 \emph{By the convexity of $\Omega$, for any $x_0 \in \Omega$,}
\begin{equation*}
\emph{AC} (\Omega)=\emph{AC}_{x_0} (\Omega): =\{u\in \mathbb R^n\,|\, x_{0}+ tu\in \Omega, \mbox{ for all }t\geq 0\}.
\end{equation*}
\emph{Note that $\mbox{AC} (\Omega )$ is a properly convex closed
cone in $\mathbb R^n$ if $\Omega $ is properly convex.}
\begin{remark} \label{rem-aff}
\begin{enumerate}
\item [\rm (i)]{It is well known that there is no bounded quasi-homogeneous affine domain (see
[1]). So we see that $\mbox{AC} (\Omega )$ is not empty if
$\Omega$ is quasi-homogeneous.}
\item [\rm (ii)]{Actually, the asymptotic cone of $\Omega $ is the maximal closed cone which can be contained in $\Omega$.}
\end{enumerate}
\end{remark}
\end{definition}

Now we will show that any properly convex quasi-homogeneous affine
domain is foliated by cosets of its asymptotic cone. To prove
this, we need the following result (see \cite{V3}).

\begin{theorem}[Vey]\label{thm-vey}
Let $\Omega$ be a properly convex open cone in $\mathbb R^n$ and
$G<\mbox{GL}(n,\mathbb R)$ acts syndetically on $\Omega$. If $L$
is a $G$-stable subspace generated by the cone $\overline{\Omega
}\cap L$, the sections $\Omega_{x}=\Omega \cap (x+L)$ are cones
for each $x\in \Omega$.
\end{theorem}

Let $L$ be the linear subspace of $\mathbb R^n$ which is generated
by
 $\mbox{AC} (\Omega )$ and $\Omega_{x}$ the intersection of $\Omega$ and the affine subspace $x+L$, that is,
 $\Omega_{x}=\Omega \cap x+L$.
 Then we see immediately that $x+\mbox{AC} (\Omega ) \subset \Omega_{x}$ for all $x \in \Omega$ since $\Omega$ is convex.
 Actually we can show that $\Omega_{x}$ is the translation of $\mbox{AC} (\Omega )$ for all $x \in
 \Omega$.
 More precisely, any section of $\Omega$ which is cut in parallel with $L$ is itself a cone, which is exactly an asymptotic cone.
\begin{theorem}\label{thm-folliation}
Let $\Omega$ be a properly convex quasi-homogeneous affine domain
in $\mathbb R^n$. Then $\Omega $ admits a parallel foliation by
cosets of the asymptotic cone of $\Omega$, i.e., $\Omega_{x}$ is a
translation of the asymptotic cone of $\Omega $ for all $x\in
\Omega$.
\end{theorem}
{\bf Proof. }
 In the vector
space $\mathbb R^{n+1}=\{ (x,t)\,|\,x\in \mathbb R^n ,t\in \mathbb
R\}$, the affine hyperplane $L_{1}=\{ (x,t)\,|\,x\in \mathbb R^n
,t=1 \}$ is isomorphic to $\mathbb R^n$. So we can consider
$\Omega $ as an open set of $L_{1}$, and $G$ as a subgroup of
$\mbox{GL} (n+1,\mathbb R)$ conserving $L_1$ and $L_0$ by the
following correspondence:
\begin{center}
$(i,\rho) :\, (\mathbb R^n, \aff{n}) \to (\mathbb R^{n+1},
\mbox{GL}(n+1, \mathbb R)),$
\end{center}
\begin{align*}i(x_1, \dots,x_n) &= (x_1,  \dots,x_n,1 ),\\
                 \rho(A,a) &=\begin{pmatrix}
  A & a \\
  0 & 1
\end{pmatrix}.\
\end{align*}
We define a properly convex open cone $\Omega_1$ of $\mathbb
R^{n+1}$:
\begin{equation*}\Omega_{1}=\{(x,t)\,|\, t>0, \frac{x}{t}\in
\Omega \}.
\end{equation*}
We easily see that $\overline{\Omega_{1}}\cap L_{0}$ is isomorphic
to the asymptotic cone of $\Omega $. Let $G_{1}$ be the subgroup
of $\mbox{GL}(n+1,\mathbb R)$ generated by $G$ and $2^m
I_{n+1}\,(m\in \mathbb Z)$. Then $G_{1}$ acts on $\Omega_{1}$ and
preserve $L_{0}$. And we have the following equalities:
\begin{align*}
\Omega_{x}&=\Omega_1 \cap (x+\langle L_{0}\cap \overline{\Omega_{1}}\rangle)\\
&=\Omega \cap (x+\langle L_{0}\cap \overline{\Omega_{1}}\rangle)\\
&=\Omega \cap (x+\langle \mbox{AC}(\Omega)\rangle) \qquad \mbox{ for any $x \in \Omega$}.
 \end{align*}
 By applying Theorem \ref{thm-vey} to
$L=\langle L_{0}\cap \overline{\Omega_{1}}\rangle$, we see
$\Omega_{x}$ must be a cone and thus must be the asymptotic cone
of $\Omega$ by convexity.
 $\qquad \square $

From Theorem \ref{thm-folliation} we get the following corollary
which was originally proved by Vey \cite{V3}.
\begin{corollary}\label{coro-cone}Let $\Omega$ be a properly convex quasi-homogeneous affine domain in $\mathbb R^n$.
If the dimension of $\mbox{AC}(\Omega)$ is equal to $n$, then
$\Omega$ is a cone.
\end{corollary}
Theorem \ref{thm-folliation} implies that for each $x \in \Omega$ there exists a point $s(x)$ in its boundary such that
 $\Omega_x = \mbox{AC}(\Omega )+s(x)$.
So we get a continuous map $s\,:\,\Omega \to \partial \Omega$ and
an one parameter group of homeomorphisms of $\Omega$ by the
following formula:
\begin{equation*}
c_{t}(x)=s(x) +e^{t}(x-s(x)) \mbox{ for }t\in \mathbb R,\, x\in
\Omega.
\end{equation*}
This flow preserves every asymptotic section $\Omega_x$ and gives a nonvanishing vector field on $\Omega$ which is similar to the Euler vector field of any cone.
\begin{definition}

\begin{enumerate}
\item [\rm (i)]
\emph{We call $p\in \partial \Omega$ an } asymptotic cone point
\emph{of a properly convex affine domain $\Omega$ if $p$ is a cone
point of $\Omega_x$ for some $x \in \Omega$, that is, $p=s(x)$ for
some $x\in \Omega$.}

\item [\rm (ii)] \emph{We call $p\in \partial \Omega$ an} extreme point \emph{of
a properly convex projective domain $\Omega$ if $p$ is not
expressed by a convex sum of any two points in $\partial \Omega$,
or equivalently, a zero dimensional face of $\Omega$.}
\end{enumerate}
\end{definition}
It is well-known that if $\Omega$ is a properly convex domain in
$\mathbb {RP}^n$ then $\overline \Omega$ is the closed convex hull
of the set of its extreme points (the Krein-Milman theorem; see
\cite{Ru}). The following proposition says that asymptotic cone
points of any properly convex affine domain are the only extreme
points in $\mathbb R^n$, that is, all other extreme points are at
infinity.

\begin{proposition}Let $\Omega$ be a quasi-homogeneous properly convex
affine domain in $\mathbb R^n$. Let $E$ be the set of all extreme points of $\overline \Omega \subset \mathbb {RP}^n$ and $S$ the
set of all asymptotic cone points $s(x)$. Then $$E=S\cup E_{\infty}
$$
where $E_{\infty}$ denote the set of all extreme points of
$\overline{\mbox{AC}(\Omega)}\cap \mathbb {RP}^{n-1}_{\infty}=\overline \Omega \cap \mathbb {RP}^{n-1}_{\infty}$.
\end{proposition}
{\bf Proof. }
\begin{enumerate}
\item[\rm(i)] $S \subset E$

Suppose not. Then there exists an open line segment $l \subset \overline \Omega$ which contains
a point $s(x)$ for some $x\in \Omega$. Note that
$\overline{ab}$ is obviously contained in $\partial \Omega$ since
$s(x) \in \partial \Omega$. Let $F$ be the face containing $s(x)$.
By Lemma \ref{lem-bounded}, we know that $F$ cannot be bounded and so it
contains an asymptotic direction $v$. This means that $s(x)+tv
\in F$ for all $t\geq 0$. Since $F$ is open in $\langle F
\rangle$, we can choose a negative real number $t$ such that $s(x)+tv \in F$. But this
contradicts that $s(x)$ is an asymptotic cone point. Therefore we
conclude that $S\subset E$.
\item[\rm(ii)] $E_{\infty} \subset E$

This follows from the fact that any
extreme point of a face of $\Omega$ is an extreme point of
$\Omega$.
\item[\rm(iii)] $E \subset S \cup E_{\infty}$

Let $e\in \mathbb R^n$ be an extreme point of $\Omega$, that is,
$e\in E \setminus E_{\infty}$. Consider the affine cone
$e+\mbox{AC}(\Omega)$. The convexity of $\overline \Omega$ implies
that $e+\mbox{AC}(\Omega)$ is a subset of $\overline \Omega$. If
$(e+\mbox{AC}(\Omega)) \cap \Omega \neq\varnothing$, then
$e+\mbox{AC}(\Omega)$ must be a section of $\Omega$ by Theorem
\ref{thm-folliation} and thus $e \in S$. If $e+\mbox{AC}(\Omega)
\subset
\partial \Omega $, then $e+\mbox{AC}(\Omega)$ must be a closed face of
$\Omega$ by Theorem~\ref{thm-folliation} again. Choose a sequence
$\{ g_n\} \subset \mbox{Aut}_{\text{aff}}(\Omega)$ converging to $g \in \mbox{PM}(n+1,\mathbb R^n)$
such that
$R(g)=\{ e \}$. Then $K(g)=\partial \mathbb R^n$ by Lemma
\ref{lem-kernel} and so $g_n(s)$ converges to $e$ for any $s \in
S$. Since $g_n(s) \in S$, we get $e\in \overline {S}\cap \mathbb
R^n$.

 Suppose that $S'=(\overline {S}\cap \mathbb R^n) \setminus S$ is not empty. Since
$\overline {S}$, $\mathbb R^n$ and $S$ are preserved by $\mbox{Aut}_{\text{aff}}(\Omega)$,
so is $S'$.
Choose $s_0 \in S $ and a sequence $\{ g_n\} \subset
\mbox{Aut}_{\text{aff}}(\Omega)$ converging to $g \in
\overline{\mbox{Aut}(\Omega)} < \mbox{PM}(n+1,\mathbb R^n)$ such
that $R(g)=s_0$. Then $K(g)=\partial \mathbb R^n$ and any point of
$\mathbb R^n$ converges to $s_0$. But all point of $S'$ cannot
converge to $s_0$ since $g_n( S')=S'$ for all $g_n$ and $S$ is an
open subset of $\overline S$. Therefore we conclude $e \in S
$.$\qquad \square $
\end{enumerate}

\section{Strictly convex quasi-homogeneous domains}\label{SC}
We will show in this section that there is only one strictly
convex quasi-homogeneous affine domain up to affine equivalence,
although there are infinitely many quasi-homogeneous strictly
convex domains in $\rp{n}$ as mentioned in the Section \ref{pre}.
To do this, we first investigate strictly convex quasi-homogeneous
projective domains since every strictly convex quasi-homogeneous
affine domain is a strictly convex quasi-homogeneous projective
domain, as will be shown later in this section. For example, a
strictly convex affine domain $D=\{ (x,y) \in \mathbb R^2
\,|\,y>1/x\}$ is not strictly convex when it is considered as a
projective domain, and it is easy to show that $D$ is not
quasi-homogeneous.

\subsection{Strictly convex quasi-homogeneous projective domains with twice differentiable boundary}\label{strict1}
The aim of this section is to prove
\begin{theorem}\label{thm-ball}
Let $\Omega$ be a strictly convex quasi-homogeneous projective
domain with twice differentiable boundary. Then $\Omega$ is
projectively equivalent to a ball.
\end{theorem}
This was already proved in two dimensional case by Kuiper
\cite{Ku} and Benz\'{e}cri \cite{B}. For arbitrary dimensional
cases, we need some lemmas.
\begin{lemma}
Let $f$ and $g$ be twice differentiable real-valued functions
defined on $(-r, r)$. Suppose $f(0)=g(0)$, $f'(0)=g'(0)$ and
$f(x)\geq g(x)$ for all $x \in $ $(-r, r)$. Then $f''(0) \geq
g''(0)$.
\end{lemma}
{\bf Proof. } If $f'(x) \leq g'(x)$ on $[0, a)$ for some $a > 0$,
then $f-g$ is a non-increasing function on $[0, a)$. Since $f-g
\geq 0$ and $(f-g)(0)=0$ by hypothesis, $(f-g)(x)=0$ on $[0,a)$
and this implies $f''(0)=g''(0)$.

Otherwise, for each integer $n > 0$ there exists a positive real
number $p_n$ such that $p_n < 1/n$ and $f'(p_n) > g'(p_n)$. Twice
differentiability of $f-g$ implies the following equalities.
$$f''(0)-g''(0) = \lim_{h \to 0} \frac{f'(h)-g'(h)}{h}=\lim_{n \to
\infty} \frac{f'(p_n)-g'(p_n)}{p_n}$$ Since for any
$n$$$\frac{f'(p_n)-g'(p_n)}{p_n}
> 0 ,$$   we conclude $f''(0)-g''(0) \geq 0$.
 $\qquad \square $

\begin{lemma}\label{lem-hessian}
Let $f$ be a strictly convex twice differentiable function defined
on an open subset $D$ of $\mathbb R^n$. Then there exists a point
$x \in D$ such that the Hessian of $f$ at $x$ is positive
definite. \footnotemark  \footnotetext{The author learned the idea
of the proof from Prof. Mohammad Ghomi.}
\end{lemma}
{\bf Proof. } We can choose a hemisphere which lies below and
touch the graph of $f$ at some $x \in D$. Let $g$ be the defining
function of the hemisphere. Then the Hessian of $g$ is positive
definite everywhere. By the above lemma, we see $D_v^2f(x) \geq
D_v^2g(x) > 0 $ for any direction $v$ and this yields the desired
result.
 $\qquad \square $

{\bf The proof of Theorem \ref{thm-ball}}:  First we choose an
affine chart containing $\Omega$. Let $f$ be a local boundary
equation of $\Omega$. By Lemma~\ref{lem-hessian}, there exists a
point $p$ in $\partial \Omega$ such that the Hessian of $f$ at $p$
is positive definite. By Corollary 6 of Rockafellar \cite{R3}, any
twice differentiable convex function has a second-order Taylor
expansion and thus we can take a suitable basis of $\mathbb R^n$
such that $p=(0,\dots,0)$ and
\begin{equation*}
\lim_{(x_{1},\dots,x_{n-1})\to 0}
\frac{f(x_{1},\dots,x_{n-1})}{{x_1}^2+ \dots +{x_{n-1}}^2}=1.
\end{equation*}
This means that $\Omega$ has an osculating ellipsoid and so
$\Omega $ is projectively equivalent to a ball by Theorems
\ref{thm-benz} and \ref{prop-benz} (ii). $\qquad \square $

The above proof says in fact the following.
\begin{proposition}\label{loc-ball}
Let $\Omega$ be a strictly convex quasi-homogeneous projective
domain whose boundary is twice differentiable on some open subset.
Then $\Omega$ is projectively equivalent to a ball.
\end{proposition}
\subsection{Strictly convex quasi-homogeneous affine domains}\label{strict2}
In this section we will show that an n-dimensional paraboloid is
the only strictly convex affine domain in $\mathbb R^n$ up to
affine equivalence, as stated in the beginning of this section
\ref{SC}. We will obtain this result through the following steps.

\begin{enumerate}
\item[\rm(I)]Every strictly convex quasi-homogeneous affine domain has an 1-dimensional asymptotic
cone, that is, it has an invariant direction.
\item[\rm(II)] The boundary of any strictly convex quasi-homogeneous affine domain is
 homogeneous and thus of class $C^{\infty}$.
\item[\rm(III)] A strictly convex quasi-homogeneous affine domain is affinely equivalent to a paraboloid.
\end{enumerate}

\begin{proposition}\label{prop-strictly}
Let $\Omega $ be a  quasi-homogeneous properly convex affine
domain in $\mathbb R^n$. Then $\Omega$ is strictly convex if and
only if $\mbox{AC} (\Omega )$ is one-dimensional.
\end{proposition}
{\bf Proof. } Suppose $\Omega $ is strictly convex and $\dim
\mbox{AC}(\Omega )\geq 2$. Consider $\Omega$ as a projective
domain. Then $\overline{\mbox{AC}(\Omega )}\cap \rp{n-1}_{\infty}$ is a closed face
of $\Omega $ with $\dim \geq 1$, and thus $\Omega \subset \rp{n}$
has a line segment $\ell$ in $\partial_p \Omega$.
(Note that this implies that $\Omega$ is not strictly convex as a projective domain.)
 From Theorem \ref{prop-benz} (iv),
 we see $\Omega$ has a triangular section, and thus $\partial_a \Omega \cap \mathbb
R^n$ has also a line segment. But this contradicts the fact that
$\Omega$ is a strictly convex affine domain.

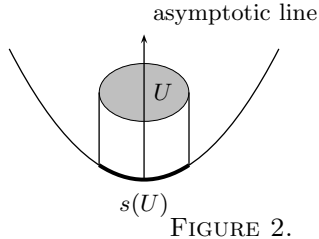
\begin{figure}[ht]
\begin{center}
 \begin{pspicture}(-4.5,0.1 )(7,2.6)
      \psset{xunit=17pt,yunit=22pt}
      \parabola[linewidth=0.5pt](3,2.25)(0,0)
      \parabola[linewidth=1.5pt](1,0.25)(0,0)
      \psellipse[linestyle=dotted,dotsep=1.5pt](0,1.5)(1,0.5)
      \psellipse[linewidth=0.1pt,
      fillstyle=solid,fillcolor=lightgray](0,1.5)(1,0.5)
      \psline[linewidth=0.5pt](1,1.5)(1,0.25)
      \psline[linewidth=0.5pt](-1,1.5)(-1,0.25)
      \psline[linewidth=0.5pt]{->}(0,0)(0,2.5)
      \uput[45](0,1.2){\small{$U$}}
      \uput[270](0,-0,5){\small{$s(U)$}}
      \uput[45](0,2.5){\small{asymptotic line}}
    \end{pspicture}
    \caption{} \label{firstfigure}
\end{center}
\end{figure}

To prove the other direction, suppose that $\dim(\mbox{AC}
(\Omega))=1$ and $\Omega$ is not strictly convex. Suppose that
there exists an open neighborhood $U \subset \mathbb R^n$ of $p
\in \partial_a \Omega$ such that $\partial_a \Omega \cap U$ has no
line segment. Then $p$ is an extreme point of $\Omega$ and there
exist $\{g_i\} \subset \mbox{Aut}(\Omega)$ such that $\lim_{i \to
\infty} g_i =g$ and $R(g) = \{p\}$ by Lemma~\ref{lem-saillant}.
Note that $\overline{\Omega} \cap \partial \mathbb R^n$ is one
point set which equals to $\overline{\mbox{AC}(\Omega)} \cap
\partial \mathbb R^n$. For any open subset $V$ of  $\mathbb R^n$,
$g_i(V)$ converges uniformly to $p$, since $\mathbb R^n$ does not intersect $K(g)$ by Lemma
\ref{lem-kernel}. (In fact, $K(g)$ is exactly
$\partial \mathbb R^n = \rp{n-1}_{\infty}$.)
 Thus $g_i(V) \subset U$ for
sufficiently large $i$ and this implies $g_i(V) \cap \partial_a
\Omega$ has no line segment and so does $V \cap
\partial_a \Omega$ if it is not empty. Since $V$ is arbitrary we
conclude that $\partial_a \Omega$ has no line segment, which is a
contradiction. So $\Omega$ can not be strictly
convex anywhere, i.e., any affine boundary point of $\Omega$ does not have a
neighborhood containing no line segment.

Let $L$ be an one dimensional linear subspace generated by
$\mbox{AC}(\Omega)$. For each $x \in \Omega$, we defined $s(x)$ in
section \ref{AC} by the point in $\partial_a \Omega$ satisfying
the following: $$s(x) + \mbox{AC}(\Omega) = \Omega \cap (x+ L).$$
Recall that the map $s\,:\, \Omega \to
\partial_a \Omega $ is continuous.  Now choose an open set $U$
contained in $\Omega$.  Then $s(U)$ is an open subset of $\partial_a
\Omega$. Since $\Omega$ is not strictly convex everywhere, $s(U)$
contains an open line segment $l$.  But by \lemref{lem-bounded},
$\partial_a \Omega$ can't have any nonzero dimensional bounded face.
Therefore the asymptotic half line must be contained in the face
containing $l$. But this contradicts that every asymptotic line
starting at a point in $s(U)$ passes $U\subset \Omega $ (see the
figure~\ref{firstfigure}).
 $\qquad \square $
\begin{remark}
This proposition implies that every strictly convex
quasi-homogeneous affine domain is also strictly convex as a
projective domain.

\end{remark}

\begin{proposition}\label{thm-homo}
Let $\Omega $ be a quasi-homogeneous affine domain in $\mathbb
R^n$. If $\Omega $ is strictly convex, then its boundary is
homogeneous.
\end{proposition}
{\bf Proof. }Let $G$ be a subgroup of $\aff{n}$ which acts
syndetically on $\Omega$.
 By \propref{prop-strictly}, there exists a one-dimensional vector space
$L$ such that $\mbox{AC} (\Omega )=\Omega \cap L$ is a half line.
For each $x\in \Omega $, the intersection $\Omega_{x}=\Omega \cap
(x+L)$ is a half line and $s(x)$ is the starting point of
$\Omega_{x}$. As stated in section \ref{AC}, the
map $s\,:\, \Omega \to
\partial_a \Omega $ is continuous and induces a one-parameter group
of homeomorphism of $\Omega$ by the following formula :
\begin{equation*}
c_{t}(x)=s(x) +e^{t}(x-s(x)) \mbox{ for }t\in \mathbb R,\, x\in
\Omega.
\end{equation*}
It is easy to show that $c_{t}$ is $\mbox{d}_{H}$-decreasing for
$t\geq 0$ (see figure \ref{secondfigure}) when $\Omega$ is
strictly convex, that is,
\begin{equation*}
\mbox{d}_{H}\bigl(c_{t}(x),c_{t}(y)\bigr)\leq \mbox{d}_{H}(x,y)
\mbox{ for all } t\geq 0.
\end{equation*}

\begin{figure}[ht]
\begin{center}
 \begin{pspicture}(-4.5,0.6 )(7,3)
      \psset{xunit=22pt,yunit=17pt}
      \parabola[linewidth=0.5pt](3,4.5)(0,0)
      \psline[linewidth=0.5pt](6,4)(-3,-0.5)
      \psline[linewidth=0.5pt](6,4)(-3,1.75)
      \psline[linewidth=0.5pt](6,4)(-3,4)
      \psline[linewidth=0.5pt](2,2)(2,4.5)
      \psline[linewidth=0.5pt](-1,0.5)(-1,4.5)
      \psline[linestyle=dotted,dotsep=1pt](2.5,3.125)(2.5,4.25)
      \psline[linestyle=dotted,dotsep=1pt](-2,2)(-2,4.25)
      \uput[45](-0.9,1){\small{$s(x)$}}
      \uput[45](-0.9,2.4){\small{$x$}}
      \uput[45](-0.9,3.8){\small{$c_{t}(x)$}}
      \uput[135](1.9,1.8){\small{$s(y)$}}
      \uput[135](1.8,2.9){\small{$y$}}
      \uput[135](1.9,3.8){\small{$c_{t}(y)$}}
    \end{pspicture}
    \caption{} \label{secondfigure}
\end{center}
\end{figure}

On the other hand, we see that $g(\Omega_{x})=\Omega_{gx}$ since
$L$ is an invariant direction, that is, since $L$ is invariant
under the action of  the linear parts of $\text{Aut}_{\text{aff}}
( \Omega )$. This implies
\begin{equation}\label{eq-sg}
s(gx)=gs(x) \mbox{   for all } g\in G.
\end{equation}
 We may assume $G$ is closed in $\mbox{Aff}(n,{\mathbb R})$. Then
 $G$ acts on $\Omega$ properly by Proposition \ref{prop-acts}. Now we can define $c_{t}$ on the orbit
space $G\setminus \Omega $ by \eqref{eq-sg}. In fact, for all
$g\in G$ and $y\in \Omega$
\begin{equation*}
\begin{split}
c_t(gy)&=s(gy)+e^{t}(gy-s(gy))=gs(y)+e^t(gy-gs(y))\\
&=e^tgy+(1-e^t)gs(y)=g(e^ty+(1-e^t)s(y))= gc_t(y)
\end{split}
\end{equation*}
  since every
affine transformation preserves a convex combination.
 Note that
\lemref{lem-exist} implies that there exists $g_{x,y}\in G$ such
that
\begin{equation}\label{eq-metric2}
\mbox{d}_{G}(Gx,Gy)=\mbox{d}_{H}(x,g_{x,y}y)
\end{equation}
where $\mbox{d}_{G}$ is defined in \eqref{eq-metric}.

It is obvious that for $t\geq 0$, the homeomorphisms $c_{t}$ of
$G\setminus \Omega$ are distance decreasing. By the way, the fact
that $G$ acts on $\Omega$ syndetically implies that $G\setminus
\Omega $ is a compact metric space and it is well known that every
distance decreasing surjection from a compact metric space onto
itself is an isometry.
 But $c_{t}$ is
strictly distance decreasing on $\Omega $ for any pair $(x,y)$
such that $s(x)\neq s(y)$, that is,
$\mbox{d}_{H}(x,y)>\mbox{d}_{H}\bigl(c_{t}(x), c_{t}(y)\bigr)$.

 Suppose there
exists a pair $(x,y)$ such that $s(x)\neq gs(y)$ for all $g \in
G$. Then we have
\begin{equation*}
\mbox{d}_{G}(Gx,Gy)=\mbox{d}_{H}(x,g_{x,y}y)>\mbox{d}_{H}\bigl(c_t(x),c_t(g_{x,y}y)\bigr)=\mbox{d}_{H}\bigl(
c_t(x),g_{x, y}c_t(y)\bigr)
\end{equation*}
and this implies
\begin{equation*}
 \mbox{d}_{G}(Gx,Gy)
>\mbox{d}_{G} \bigl( c_{t}(Gx),c_{t}(Gy)\bigr).
\end{equation*} This is a
contradiction. Therefore, for any pair $(x,y)$ we get $s(x)=gs(y)$
for some $g \in G$ and so we conclude that $\partial_a \Omega $ is
homogeneous.
 $\qquad \square $
\begin{remark}
We showed in the proof of the above proposition that the set of
all asymptotic cone points $S$ is homogeneous when it is strictly
convex. In this case, $S=\partial \Omega$ since $\Omega$ has a
1-dimensional asymptotic cone. In case $\dim \emph{AC} (\Omega )$
is not 1-dimensional, $S$ is a proper subset of $\partial \Omega$.
For example, if $\Omega=\{(x,y,z) \in \mathbb R^3\,|\,x>0,
y>z^2\}$ then $S=\{(0,y,z) \in \mathbb R^3\,|\, y=z^2\}$. $S$
seems to be still homogeneous even in the case of $\dim \emph{AC}
(\Omega ) > 1$, which implies that the Lie group $\emph{Aut}
(\Omega)$ is at least of dimension $n-\dim \emph{AC} (\Omega )$.
But this is not proved yet.
\end{remark}

Now, we get the following main theorem of this section.
\begin{theorem} \label{thm-para}Every strictly convex quasi-homogeneous affine domain in
$\mathbb R^n$ is affinely equivalent to an n-dimensional
paraboloid.
\end{theorem}
{\bf Proof. } We have shown in Theorem \ref{thm-homo} that
$\partial_a \Omega$ is homogeneous and we can see $\partial_a
\Omega$ is twice differentiable almost everywhere by Corollary 6
of \cite{R3} (in fact, $\partial_a \Omega$ is $C^{\infty}$ since
it is the orbit of the Lie group Aut($\Omega$)). But the
homogeneity of $\partial_a \Omega$ implies that $\partial_a
\Omega$ is twice differentiable everywhere and thus $\Omega$ is
projectively equivalent to a ball in $\rp{n}$ by Proposition
\ref{loc-ball} since $\Omega$ is also strictly convex when it is
considered as a projective domain. Therefore we conclude that
$\Omega$ is affinely equivalent to a paraboloid since $\overline
{\Omega} \cap
\partial\mathbb R^n=\overline{\mbox{AC}(\Omega)} \cap
\partial\mathbb R^n$
 is a point by Proposition \ref{prop-strictly}.
 $\qquad \square $

\subsection{Strictly convex quasi-homogeneous projective domains}\label{strict3}
Using the results of the previous section, we will show that if
$\Omega$ is a strictly convex quasi-homogeneous projective domain,
then
\begin{enumerate}
\item [\rm (i)]$\partial \Omega$ is at least
$C^1$,
\item [\rm (ii)]$\Omega$ is an ellipsoid if $\partial \Omega$ is
twice differentiabe,
\item [\rm (iii)]if $\partial \Omega$ is $C^{\alpha}$ on an open subset of $\partial \Omega$, then
$\partial \Omega$ is $C^{\alpha}$ everywhere.
\end{enumerate}

 Note that we already proved (ii) in Theorem \ref{thm-ball}. Now we
 prove (i).
\begin{proposition}\label{prop-C^1}
Let $\Omega $ be a quasi-homogeneous properly convex domain in
$\rp{n}$. Then $\Omega $ is strictly convex if and only if the
boundary of $\Omega $ is continuously differentiable.
\end{proposition}
{\bf Proof. } Note that for a convex function, the
differentiability implies the continuity of the derivative (see
Corollary 25.5.1 of \cite{R1}) and the existence of all partial
derivatives implies the differentiability (see \cite{Sc} p.27).
So, to prove that $\partial \Omega$ is $C^1$ it suffices to show
that $\partial \Omega$ has partial derivatives everywhere. Suppose
that $\Omega $ is strictly convex and there exists a point $x\in
\partial \Omega$ at which $\partial \Omega$ does not have some partial derivative.
Then we can choose a 2-dimensional subspace $H$ such that $H \cap
\Omega$ is a 2-dimensional section of $\Omega$ containing $x$ in
its boundary and $\partial (H \cap \Omega)$ is not differentiable
at $x$. This means that $x$ is a conic face of $H \cap \Omega $,
hence $\Omega $ has a triangular section by Theorem
\ref{prop-benz} (iv). But this contradicts the strict convexity of
$\Omega$.

Now we prove the converse. Suppose that $ \partial \Omega$ is
differentiable and $\Omega$ is not strictly convex. Since $\Omega$
is not strict convex, $\Omega$ has a triangular section by
Corollary \ref{cor-con}. This implies that $ \partial \Omega$ can
not be differentiable at any vertex point of the triangle. This
contradiction completes the proof.
 $\qquad \square $
\begin{proposition}\label{prop-euclidean}
Let $\Omega $ be a quasi-homogeneous strictly convex domain in
$\rp{n}$. If $\partial  \Omega$ is $C^{\alpha}$ on some
neighborhood $V$ of $p\in \partial \Omega$, then $\partial \Omega$
is $C^{\alpha}$.
\end{proposition}
{\bf Proof. }
  We showed in Lemma \ref{lem-saillant} that for each point
   $p\in \partial \Omega$, there exists
 a sequence $\{g_{n}\}\subset \mbox{Aut} ( \Omega ) $ and $x\in \Omega$
  such that $g_{n}x$ converges to $p$.

Obviously $g$ is singular. By ~\lemref{lem-Ben}, $K(g)\cap \Omega
=\emptyset$ and $R(g)\cap \Omega =\emptyset$. Strict convexity of
$\Omega$ implies that the rank of $g$ is $1$ and $R(g)=\{p\}$ and
$K(g)\cap
\partial \Omega$ must be a single point, say $q$, by Lemma \ref{lem-Ben} again.
 For any $y\in \partial \Omega
\setminus \{q\}$, $g(y)=p$ and thus $g_{i}(y)\in V$. Then
$g_{i}^{-1}(V)$ is a $C^{\alpha}$neighborhood of $y$. Therefore
 $\partial \Omega$ is $C^{\alpha}$ except
$q= K(g) \cap \overline\Omega$. Suppose that $\partial \Omega$ is
not $C^{\alpha}$ at $q$. Then $q$ must be a fixed point. Since
$\partial \Omega$ is always $C^1$ by Proposition ~\ref{prop-C^1},
there exists a unique supporting hyperplane $H$ and thus $H$ is
invariant under $\mbox{Aut}( \Omega )$. So $\Omega$ must be a
strictly convex quasi-homogeneous affine domain and so it is
projectively equivalent to a ball by Theorem \ref{thm-para}. This
contradicts that $\partial \Omega$ is not $C^{\alpha}$ at $q$.
Therefore $\partial \Omega$ is $C^{\alpha}$ everywhere.
 $\qquad \square $

From Propositions ~\ref{prop-C^1} and ~\ref{prop-euclidean} and
Theorem \ref{thm-ball}, we get immediately the following
proposition.
\begin{proposition}\label{thm-ball2}
Let $\Omega$ be a strictly convex quasi-homogeneous domain in $
\rp{n}$. Then either  (i) $\Omega$ is an ellipsoid, or (ii)
$\partial \Omega$ is $C^1$ and fails to be twice differentiable on
a dense subset.

\end{proposition}

A compact convex projective manifold $M$ is a quotient of a convex
projective domain $\Omega$ in $\rp{n}$ by a discrete subgroup
$\Gamma$ of $\mbox{PGL}(n+1,\mathbb R)$ acting on $\Omega$
properly discontinuously and freely. From the above proposition we
get
\begin{corollary}
If a strictly convex projective domain $\Omega$ covers a compact
projectively flat manifold, then either $\Omega$ is an ellipsoid
or $\partial \Omega$ fails to be twice differentiable on a dense
subset.
\end{corollary}

\begin{remark}
This corollary was proved in the 2-dimensional case by Kuiper
\cite{Ku}.
\end{remark}

The following proposition shows that there is no quasi-homogeneous
strictly convex projective domain which is not divisible. In the
2-dimensional case, there is a more elementary and direct proof,
which will be given in the next section. In fact, we will show
later that every 3-dimensional quasi-homogeneous properly convex
affine cone is divisible. It seems to be true that any irreducible
properly convex projective domain is either homogeneous or
divisible. Y. Benoist \cite{Bn2}
 proved this when the domain has an irreducible automorphism group.
\begin{proposition}\label{str-aut}
Let $\Omega$ be a strictly convex quasi-homogeneous projective
domain in $ \rp{n}$. Suppose that $\Omega$ is not an ellipsoid.
Then $\mbox{Aut}(\Omega$) is discrete and so $\Omega$ is a
divisible domain.
\end{proposition}
{\bf Proof. } By the Proposition 4.2 of \cite{Bn2}, it suffice to
show that Aut($\Omega$) is irreducible. Suppose Aut($\Omega$) is
reducible and $L$ is a projective subspace of $ \rp{n}$ which is
invariant under the action of Aut($\Omega$). Then $L \cap \Omega =
\emptyset$ and $L \cap \overline \Omega \neq \emptyset$, since
$\Omega$ is quasi-homogeneous So the strict convexity of $\Omega$
implies that $L \cap \Omega$ must be a point $\{ p \}$ in
$\partial \Omega$.
 This implies $p$ is a fixed point of Aut($\Omega$).
Since every strictly convex quasi-homogeneous projective domain is
of class $C^1$, $p$ has a unique tangent plane $H_p$ and it is
preserved by Aut($\Omega$). This means that $\Omega$ is a
quasi-homogeneous strictly convex affine domain and so it must be
an ellipsoid by Theorem \ref{thm-ball2}, which gives a
contradiction. $\qquad \square $

This Proposition shows that every strictly convex
quasi-homogeneous domain covers a compact projectively flat
manifold. In addition, W. Goldman and S. Choi \cite{SG} proved
that any strictly convex $ \rp{2}$-structure on a closed surface is
a projective deformation of a hyperbolic structure. A
3-dimensional case was shown by I. Kim \cite{inkang} and then Y.
Benoist \cite{Bn3} has proved it
 for more general case.

\section{Quasi-homogeneous convex  affine domains of dimension $\leq 3$}

In this section, we will characterize low dimensional
quasi-homogeneous convex affine domains which are not cones.
\begin{proposition}\label{thm-dim2}Let $\Omega$ be a properly convex
quasi-homogeneous affine domain in $\mathbb R^2$. Then $\Omega$ is
affinely equivalent to either a quadrant or a parabola.
\end{proposition}
{\bf Proof. } If dimAC($\Omega)=1$ then $\Omega$ is affinely
equivalent to a parabola by Proposition \ref{prop-strictly} and
Theorem \ref{thm-para}. Otherwise dimAC($\Omega)=2$, and hence
$\Omega$ must be a properly convex cone by Corollary
\ref{coro-cone}. Any two dimensional properly convex cone is
affinely equivalent to a quadrant, which completes the proof.
 $\qquad \square $

\begin{proposition}\label{thm-dim3}Let $\Omega$ be a properly convex quasi-homogeneous
 affine domain in $\mathbb R^3$. Suppose $\Omega $ is
not a cone. Then $\Omega $ is affinely equivalent to one of the
following:
\begin{enumerate}
\item[\rm(i)] A 3-dimensional paraboloid, i.e., $\{(x,y,z) \in \mathbb R^3\,|\, z>x^{2} +y^{2}
\}$.
\item[\rm(ii)] A convex domain bounded by a parabola $\times \mathbb R^{+}$,
 i.e.,
  $$\{(x,y) \in \mathbb R^2\,|\,y>x^{2}\}\times \{z \in \mathbb R\,|\, z>0\}=
       \{(x,y,z) \in \mathbb R^3\,|\, y>x^2, z>0\}.$$
\end{enumerate}
\end{proposition}
{\bf Proof. }By Theorem \ref{thm-para} and Proposition
\ref{prop-strictly}, $\Omega$ is affinely equivalent to a
paraboloid if $\dim \mbox{AC} (\Omega)=1$. Since $\Omega $ is
not a cone, $\dim \mbox{AC} (\Omega)\neq 3$ by Corollary
\ref{coro-cone}. Hence it suffices to show that $\Omega$ is
affinely equivalent to a convex domain bounded by a parabola
$\times \,\mathbb R^{+}$ if $\dim \mbox{AC} (\Omega)=2$. We know
that $\Omega $ admits a parallel foliation by cosets of the
asymptotic cone of $\Omega$ by Proposition \ref{thm-folliation}.
Put $\mbox{AC} _{\infty}(\Omega)=\overline{\mbox{AC} (\Omega)}\cap
\rp{2}_{\infty}
 =\overline{ab}$. Now we will show
that one of $\{a,b\}$ is a conic point of $\Omega$. Choose an
asymptotic cone point $p\in
\partial \Omega $, then the triangle $\triangle pab$ is a $2$-dimensional
section of $\Omega$. Therefore there exists supporting hyperplanes
$H_{pa}$ and $H_{pb}$ tangent to $\overline {pa}$ and
$\overline{pb}$ respectively. Let $L=H_{pa} \cap H_{pb}$. Suppose
every asymptotic cone point is contained in $L$. Then the fact that $\Omega$ is not a cone implies that the set
of all asymptotic cone points in $L\cap
\partial \Omega$ cannot be one point set and thus it is a bounded line,
 since $L$ is not an asymptotic
direction. This contradicts Lemma \ref{lem-bounded}. So, there
exists an asymptotic cone point $p'\in
\partial \Omega$ such that $p'\notin L$. This means that $p'\notin
H_{pa}$ or $p'\notin H_{pb}$. We may assume that $p'\notin
H_{pa}$. Note that $H_{pa}$ and $H_{p'a}$ are neither the same
hyperplane nor parallel. So $a$ is a conic face of $\Omega$ when
it is considered as a projective domain. In fact $a$ is the
intersection of the three independent supporting hyperplanes
$H_{pa},H_{p'a},H_{ab}$, where
 $H_{ab}$ is in fact $\partial \mathbb R^3=\rp{2}_{\infty}$.
 By Theorems \ref{thm-benz} and \ref{prop-benz} (i), we conclude that
$\Omega=a+\Omega'$ for some $2$-dimensional quasi-homogeneous
properly convex projective domain $\Omega'$. $\Omega'$ must
contain all asymptotic cone points of $\Omega$ and be an invariant
quasi-homogeneous 2-dimensional affine domain with 1-dimensional
asymptotic cone. By Theorem \ref{thm-dim2}, $\Omega'$ is a
parabola. Consequently, we have that $\Omega=\overline {pa} \times
\Omega'$.
 $\qquad \square $

If $\Omega$ is a quasi-homogeneous convex affine domain which is
not properly convex, then $\Omega$ is a product $\Omega' \times
\mathbb R^k$ with $k>0$, where $\Omega'$ is a ($n-k$)-dimensional
quasi-homogeneous properly convex affine domain. So we get the
following two corollaries from Propositions \ref{thm-dim2} and
 \ref{thm-dim3}.
\begin{corollary}Let $\Omega$ be a quasi-homogeneous convex affine domain in $\mathbb
R^2$. Then $\Omega$ is affinely equivalent to one of the
following :
\begin{enumerate}
\item[\rm(i)] $\{(x,y) \in \mathbb R^2\,|\,y>x^2\}$,
\item[\rm(ii)] $\{(x,y) \in \mathbb R^2\,|\,x>0,y>0\}$,
\item[\rm(iii)]$\{(x,y) \in \mathbb R^2\,|\,y>0\}$,
\item[\rm(iv)] $\mathbb R^2$.
\end{enumerate}
\end{corollary}

Every 3-dimensional quasi-homogeneous properly convex cone is
either a strictly convex cone or a simplex cone since a
2-dimensional quasi-homogeneous properly convex projective domain
which is not a triangle is strictly convex. So we get
\begin{corollary}\label{qh3} Let $\Omega$ be a convex quasi-homogeneous domain in $\mathbb
R^3$. The $\Omega $ is affinely equivalent to one of the
following :
\begin{enumerate}
\item[\rm(i)] $\{(x,y,z) \in \mathbb R^3\,|\,z>x^2+y^2\}$,
\item[\rm(ii)] $\{(x,y,z) \in \mathbb R^3\,|\,y>x^2,z>0\}$,
\item[\rm(iii)]$\{(x,y,z) \in \mathbb R^3\,|\,x>0,y>0\}$,
\item[\rm(iv)] $\{(x,y,z) \in \mathbb R^3\,|\,y>x^2\}$,
\item[\rm(v)] $\{(x,y,z) \in \mathbb R^3\,|\,z>0\}$,
\item[\rm(vi)] $\mathbb R^3$,
\item[\rm(vii)] a quasi-homogeneous strictly convex affine cone,
\item[\rm(viii)]a simplex cone.
\end{enumerate}
\end{corollary}

\begin{remark}\label{str-div}
\begin{enumerate}
\item [\rm(i)]Every 3-dimensional quasi-homogeneous properly convex affine cone is
divisible. The reason is as follows: Any quasi-homogeneous
properly convex 2-dimensional projective domain which is neither a
triangle nor an ellipsoid is a strictly convex domain whose
boundary is not twice differentiable on a dense subset by
Proposition \ref{thm-ball2}. If the automorphism group of a
2-dimensional strictly convex domain has a 1-dimensional Lie
subgroup, the boundary has a 1-dimensional orbit and consequently
is twice differentiable, that is, the boundary is twice
differentiable on some open set. So the automorphism group of a
2-dimensional quasi-homogeneous strictly convex domain whose
boundary is not twice differentiable, must be discrete. Therefore
every 2-dimensional quasi-homogeneous properly convex projective
domain is divisible and so is any 3-dimensional quasi-homogeneous
properly convex affine cone. Actually, this fact follows from
Proposition \ref{str-aut}.
\item [\rm(ii)]Any strictly convex cone of (vii) is either an elliptic cone
$\{(x,y,z) \in \mathbb R^3\,|\,x^2-y^2-z^2>0\}$ or its boundary is
not twice differentiable on any open subset by Proposition
\ref{thm-ball2}.
\end{enumerate}
\end{remark}
\section{Compact convex affine manifolds of dimension $\leq4$}

Let $M$ be a compact convex affine manifold. Then we get a divisible domain $\Omega$ and a discrete subgroup $\Gamma$ of
$\mbox{Aut}_{\text{aff}}(\Omega)$ such that $\Omega / \Gamma = M$. If $\Omega$ is properly convex, then it must be a cone by a
result of Vey \cite{V3}.

For a $2$-dimensional compact convex affine manifold $M$, we know
that $\Omega$ is a quadrant if it is properly convex. Otherwise
$\Omega$ is a half space since the only $1$-dimensional
quasi-homogeneous affine domain is $\mathbb R^{+}$.

Using the results in the previous section, we get the following.
\begin{theorem}\label{thm-dimension3}Let $\Omega$ be an affine
domain in $\mathbb R^3$ which covers a compact convex affine
manifold. Then $\Omega$ is affinely equivalent to one of the
following :
\begin{enumerate}
\item[\rm(i)] $\mathbb R^3$,
\item[\rm(ii)] $\mathbb R^2 \times \mathbb R^{+}$,
\item[\rm(iii)] $\mathbb R \times  \{(x,y) \in \mathbb R^2\,|\,y>x^{2}\} $,
\item[\rm(iv)] $\mathbb R \times \{(x,y) \in \mathbb R^2 \,|\,x>0, y>0\}$,
\item[\rm(v)] a quasi-homogeneous strictly convex cone,
\item[\rm(vi)] a simplex cone.
\end{enumerate}
\end{theorem}
{\bf Proof. }Consider a 3-dimensional quasi-homogeneous convex
affine domain $\Omega$. It is isomorphic to $\mathbb R^k \times
\Omega'$ for some ($3-k$)-dimensional properly convex affine
domain $\Omega'$. If $k=3$, $\Omega=\mathbb R^3$ and if $k=2$,
$\Omega=\mathbb R^2 \times \mathbb R^{+}$. Since a 2-dimensional
quasi-homogeneous properly convex domain is either a quadrant or a
parabola, $\Omega$ is either (iii) or (iv) if $k=1$. If $k=0$,
that is, if $\Omega$ is properly convex, then it is a cone and
thus either (v) or (vi) by Corollary \ref{qh3}.

All of these quasi-homogeneous domains are actually divisible. The
divisibility of (iii) was shown by W. Goldman in \cite{G4} and the
divisibility of (v) is implied by Remark \ref{str-div}. The
remaining cases are easy.
 $\qquad \square $

Similarly, we get the following theorem in 4-dimensional case.
\begin{theorem}\label{thm-dimension4}Let $\Omega$ be an affine
domain in $\mathbb R^4$ which covers a compact convex affine
manifold. Then $\Omega$ is affinely equivalent to one of the
following :
\begin{enumerate}
\item[\rm(i)] $\mathbb R^4$,
\item[\rm(ii)] $\mathbb R^3 \times \mathbb R^{+}$,
\item[\rm(iii)] $\mathbb R^2 \times  \{(x,y)\in \mathbb R^2\,|\,y>x^{2}\} $,
\item[\rm(iv)] $\mathbb R^2 \times  \{(x,y)\in \mathbb R^2\,|\,x>0, y>0\}$,
\item[\rm(v)] $\mathbb R \times \{(x,y,z)\in \mathbb R^3 \,|\,z>x^2+y^2\}$,
\item[\rm(vi)] $\mathbb R \times \mathbb R^{+} \times  \{(x,y)\in \mathbb R^2\,|\,y>x^{2}\} $,
\item[\rm(vii)] $\mathbb R \times $ a divisible properly convex cone of dimension $3$,
\item[\rm(viii)] a divisible
properly convex cone of dimension $4$
\end{enumerate}
\end{theorem}

W. Goldman and M. Hirsch proved in \cite{GH} that any compact
complete affine 3-manifold has parallel volume (see Cor. 3.3).
They also proved that if a compact affine n-manifold $M$ with
holonomy group $\Gamma$ has parallel volume, then the algebraic
hull $A(\Gamma)$ of $\Gamma$ acts transitively and so $\Gamma$
preserves no proper algebraic subset of $\mathbb R^n$. Considering
the developing image $\Omega$ of compact convex affine manifold of
dimension $\leq 4$, we see that the holonomy group of compact
convex affine manifold has a proper invariant algebraic set if its
developing map is not surjective. So we get the following
corollary.

\begin{corollary}Let $M$ be a compact convex affine n-manifold.
\begin{enumerate}
\item [\rm(i)]  For $n = 3$, $M$ is complete if and only if it has parallel volume.
\item [\rm(ii)] For $n = 4$, $M$ is complete if it has parallel volume.
\end{enumerate}
\end{corollary}
\section*{Acknowledgements }
I would like to thank Prof. Hyuk Kim who has encouraged and
advised me all the while I have been studying this topic. I am
also grateful to Prof. Inkang Kim for reading this article and
helpful discussions and Prof. Suhyoung Choi for his critical
comments. Conversations with Prof. Yves Benoist and Mohammad Ghomi
by e-mail have been also helpful. Finally I thank the referee and
Prof. Seong-deog Yang for suggesting a number of improvements on
the exposition.

\end{document}